\documentclass[11pt]{article}
\usepackage{mathrsfs}
\usepackage{}
\usepackage{amsfonts}
\usepackage{amssymb}
\usepackage{graphicx}
\usepackage{amsmath}

\def\sqr#1#2{{\vcenter{\vbox{\hrule height.#2pt
              \hbox{\vrule width.#2pt height#1pt \kern#1pt \vrule
width.#2pt}
              \hrule height.#2pt}}}}
\def\signed #1{{\unskip\nobreak\hfil\penalty50
              \hskip2em\hbox{}\nobreak\hfil#1
              \parfillskip=0pt \finalhyphendemerits=0 \par}}
\def\endpf{\signed {$\sqr69$}}
\def\dbN{{\mathbb{N}}}

\def\3n{\negthinspace \negthinspace \negthinspace }
\def\2n{\negthinspace \negthinspace }
\def\1n{\negthinspace }

\def\ms{\medskip}

\def\q{\quad}

\def\({\Big (}
\def\){\Big )}
\def\[{\Big[}
\def\]{\Big]}

\def\be{\begin{equation}}
\def\bel{\begin{equation}\label}
\def\ee{\end{equation}}
\def\bea{\begin{eqnarray}}
\def\eea{\end{eqnarray}}
\def\bt{\begin{theorem}}
\def\et{\end{theorem}}
\def\bc{\begin{corollary}}
\def\ec{\end{corollary}}
\def\bl{\begin{lemma}}
\def\el{\end{lemma}}
\def\bp{\begin{proposition}}
\def\ep{\end{proposition}}
\def\br{\begin{remark}}
\def\er{\end{remark}}
\def\ba{\begin{array}}
\def\ea{\end{array}}
\def\bd{\begin{definition}}
\def\ed{\end{definition}}

\newtheorem{lemma}{Lemma}[section]
\newtheorem{remark}{Remark}[section]

\newtheorem{theorem}{Theorem}[section]
\newtheorem{corollary}{Corollary}[section]

\newtheorem{definition}{Definition}[section]
\newtheorem{proposition}{Proposition}[section]

\makeatletter
   
   \@addtoreset{equation}{section}
\makeatother
\allowdisplaybreaks
\begin{document}

\title{\bf Null Controllability of Some
Degenerate Wave Equations\thanks{This work is
supported by the NSF of China under grants 11371084, 11471070 and 11171060,  by the Fundamental Research
 Funds for the Central Universities under grants 14ZZ2222 and 2412015BJ011, by the National Basic Research Program of China (973 Program) under grant 2011CB808002, and by  the Fok Ying Tong Education Foundation under grant 141001.}}

\author{Muming Zhang\thanks{School of Mathematics and Statistics, Northeast Normal
University, Changchun 130024, China. E-mail address:
zhangmm352@nenu.edu.cn. \ms } \q
 and \q Hang Gao\thanks{School of Mathematics and Statistics, Northeast Normal
University, Changchun 130024, China. E-mail address:
hanggao2013@126.com
 \ms }}

\date{}

\maketitle

\begin{abstract}
This paper is devoted to a study of  the null controllability problems
for  one-dimensional linear degenerate wave equations through a boundary controller.
First,  the well-posedness of  linear degenerate
wave equations is  discussed.  Then   the null  controllability of  some degenerate wave equations is established,
when a control acts on the non-degenerate boundary.  Different from the known controllability results in the case that a control acts on the degenerate boundary,
any initial value in state space is controllable in this case.
Also,  an explicit expression for the  controllability time is given.  Furthermore,  a counterexample  on the controllability is given for  some other degenerate wave equations.

\end{abstract}


\section{Introduction and main results}

Let $T>0$, $L>0$ and $\alpha>0$. Set $Q=(0, L)\times(0,T)$.
Consider the following linear degenerate wave equation with a boundary controller:
\begin{equation}\label{1.1}
\left\{\begin{array}{ll}
w_{tt}-(x^\alpha w_x)_x=0& (x,t)\in Q,\\[2mm]
\left\{\begin{array}{ll}
w(0,t)=0 &(0<\alpha<1)\\[2mm]
(x^\alpha w_x)(0,t)=0 &(\alpha\geq1)\\[2mm]
\end{array}\right.&t\in(0,T),\\[2mm]
w(L,t)=\theta(t) &t\in(0,T),\\[2mm]
w(x,0)=w_0(x), \ w_t(x,0)=w_1(x)&x\in(0, L),
\end{array}\right.
\end{equation}
where $\theta\in L^2(0,T)$ is the control  variable, $(w, w_t)$ is the state variable, and $(w_0, w_1)$ is
any given initial value.

In order to study the well-posedness of (\ref{1.1}),  we introduce a linear space $H_\alpha^1(0, L)$:

\noindent (1) In the case of $0<\alpha<1$,
\begin{eqnarray*}
&&H_\alpha^1(0, L)=\big{\{}u\in L^2(0, L)\ \big{|}\ u\ \text{is absolutely\
continuous\ in}\ [0,L],
\  x^\frac{\alpha}{2} u_x\in L^2(0, L)\\
&&\quad\quad\quad\quad\quad\quad\quad\quad\quad\quad\quad\quad
\text{and}\ u(0)=u(L)=0\big\}.
\end{eqnarray*}
\noindent (2) In the case of $\alpha\geq1$,
\begin{eqnarray*}
&&H_\alpha^1(0, L)=\big{\{}u\in L^2(0, L)\ \big{|}\ u\ \text{is locally\ absolutely\
continuous\ in}\ (0, L],
\ x^\frac{\alpha}{2} u_x\in L^2(0, L)\\
&&\quad\quad\quad\quad\quad\quad\quad\quad\quad\quad\quad\quad\text{and}\ u(L)=0\big\}.
\end{eqnarray*}
Then $H_\alpha^1(0, L)$ $(\alpha>0)$ is a Hilbert space with the inner product
$$(u,v)_{H_\alpha^1(0, L)}=\int_0^L(uv+x^\alpha u_xv_x)dx,\quad  \forall\ u,v\in H_\alpha^1(0, L)$$
and associated norm $\|u\|_{H_\alpha^1(0, L)}=\|u\|_{L^2(0,L)}+\|x^\frac{\alpha}{2} u_x\|_{L^2(0,L)}$.
By Hardy-poincar\'{e} inequality (see \cite[Proposition 2.1]{car}),
it is easy to check that for $\alpha\in(0,2)$, there exists a constant
$C>0$ such that for any $u\in H_\alpha^1(0, L)$ satisfies
$$\|u\|_{L^2(0,L)}\leq C\|x^\frac{\alpha}{2} u_x\|_{L^2(0,L)}.$$
Therefore, when $\alpha\in(0,2)$, $\|u\|_{H_\alpha^1(0, L)}$ is
equivalent to the norm $\|x^\frac{\alpha}{2} u_x\|_{L^2(0,L)}$.
Also, $H_\alpha^*(0, L)$ denotes the conjugate space of $H_\alpha^1(0, L)$ and $$\|v\|_{H_\alpha^*(0, L)}=\sup\limits_{\|u\|_{H_\alpha^1(0, L)}=1}\langle u,v\rangle_{H_\alpha^1(0, L),H_\alpha^*(0, L)}.$$
 Moreover, for any $0\leq \alpha<2$,
set \begin{equation}\label{zg}T_\alpha=\frac{4}{2-\alpha}L^{\frac{2-\alpha}{2}}.\end{equation}

\medskip

The degenerate wave equation (\ref{1.1}) can describe the vibration problem of an elastic string. In a neighborhood of  an  endpoint $x=0$ of this string, the  elastic  is sufficiently small or the linear density is large enough.
 First,  we  give the definition  of  solutions of the equation (\ref{1.1}) in the sense of transposition.
\begin{definition}\label{d1.2}
 For any $(w_0, w_1)\in L^2(0,L)\times H_\alpha^*(0,L)$ and $\theta\in L^2(0, T)$,
a function $w(\cdot)\in C([0,T];L^2(0,L))$  is called a solution  of
the equation $(\ref{1.1})$ in the sense of transposition,  if for any
$\xi\in L^1(0,T;L^2(0,L))$, it holds that
\begin{eqnarray}\label{1.10}
\int_Qw\xi dxdt=\langle w_1,v(0)\rangle _{H_\alpha^*(0,L),H_\alpha^1(0,L)}
-\int_0^L w_0v_t(0)dx-L^\alpha\int_0^T\theta(t)v_x(L,t)dt,
\end{eqnarray}
where $v$ is the weak solution of the following  equation:
\begin{equation}\label{1.9}
\left\{\begin{array}{ll}
v_{tt}-(x^\alpha v_x)_x=\xi & (x,t)\in Q,\\[2mm]
\left\{\begin{array}{ll}
v(0,t)=0 &(0<\alpha<1)\\[2mm]
(x^\alpha v_x)(0,t)=0 &(\alpha\geq1)\\[2mm]
\end{array}\right.&t\in(0,T),\\[2mm]
v(L,t)=0 &t\in(0,T),\\[2mm]
v(x,T)=0,\ v_t(x,T)=0 &x\in (0,L).
\end{array}\right.
\end{equation}

\end{definition}

\medskip

\noindent Then we have the following well-posedness result for the system (\ref{1.1}).
\begin{theorem}\label{t1.3}
Let $\alpha\in(0,2)$. For any $(w_0, w_1)\in L^2(0,L)\times H_\alpha^*(0,L)$ and $\theta\in L^2(0, T)$,
the system $(\ref{1.1})$ admits a unique solution
$w\in \mathscr{K}=C([0,T];L^2(0,L))\cap C^1([0,T];H_\alpha^*(0,L))$
in the sense of transposition. Moreover,
$$\|w\|_{\mathscr{K}}\leq C\big{(}\|w_0\|_{L^2(0,L)}+\|w_1\|_{H_\alpha^*(0,L)}+\|\theta\|_{L^2(0,T)}\big{)}.$$
\end{theorem}

\medskip

 The purpose of this paper is to study the null  controllability
of  the linear degenerate wave equation (\ref{1.1}), and give an expression of
the controllability time.  The system $(\ref{1.1})$ is null controllable in time $T$,
if for any initial value
$(w_0, w_1)\in L^2(0, L)\times H_\alpha^*(0, L)$,  one can find
a control $\theta\in L^2(0,T)$ such that the corresponding solution
$w$  of (\ref{1.1}) (in the sense of transposition) satisfies that $w(T)=w_t(T)=0$ in $(0, L)$.

\medskip

The  main controllability results of this paper can be stated as follows.

\begin{theorem}\label{t1}

$(1)$ Let $\alpha\in(0,2)$.
 Then for any $T>T_\alpha$,
the system $(\ref{1.1})$ is null controllable in time $T$.  For any  $T<T_\alpha$,
the system $(\ref{1.1})$ is not null controllable.

\medskip

\noindent $(2)$ Let $\alpha\in [2,+\infty)$. Then  for any $T>0$, the system $(\ref{1.1})$ is not null controllable.
\end{theorem}

\begin{remark}
Notice that when $\alpha=0$, the system $(\ref{1.1})$ is a non-degenerate linear wave equation. By the known controllability result in $\cite{zuazua}$,  the controllability time $T^*=2L$ $($for $\alpha=0)$. Letting $\alpha$ tend to zero in $(\ref{zg})$, one can find that  $\lim\limits_{\alpha\rightarrow 0} T_\alpha=T^*$.
\end{remark}

\begin{remark}
Since the linear degenerate wave equation $(\ref{1.1})$ has time-reversibility,
the null controllability of it is equivalent to the exact controllability.
\end{remark}

Up to now, there are numerous works addressing the  controllability problems of
nondegenerate parabolic and hyperbolic equations (see e.g. \cite{cla},  \cite{cui}, \cite{7},  \cite{liu}, \cite{ypf}, \cite{zhang1}
  and the
references therein).
The controllability of some degenerate parabolic equations  was studied in the last decade (see, for instance,  \cite{1},  \cite{ppj} and the references therein).
However, very little is known for the controllability of degenerate wave equations.  In \cite{Gueye},   for any $0<\alpha<1$, the null controllability of
the following
 degenerate wave equation was considered:  \begin{equation}\label{0}
\left\{\begin{array}{ll}
w_{tt}-(x^\alpha w_x)_x=0 & (x,t)\in (0,1)\times(0,T),\\[2mm]
w(0,t)=\theta(t), \ w(1,t)=0 &t\in(0,T),\\[2mm]
w(x,0)=w_0(x),\ w_t(x,0)=w_1(x)&x\in (0,1),
\end{array}\right.
\end{equation}
where  $\theta(\cdot)$ is the control variable and it acts on the degenerate boundary.  Based on the Fourier expansion method,  the null controllability  of (\ref{0}) was discussed in \cite{Gueye}.  Indeed, for any  $T>\frac{4}{2-\alpha}$ and any  initial value $(w_0,w_1)\in \mathcal{H}=H_{\alpha}^{\frac{1-2\mu}{2}}\times H_{\alpha}^{\frac{-1-2\mu}{2}}$ with $\mu=\frac{1-\alpha}{2-\alpha}$,
there exists a control $\theta\in L^2(0,T)$ such that the solution $w$ of (\ref{0})
satisfies $w(T)=w_t(T)=0$ in $(0, 1)$.  For any $s>0$,
$$H_{\alpha}^s(0,1)=\Big\{u=\sum\limits_{n\in\mathbb{N}^*}a_n\Phi_n\ \big{|}\
\|u\|_s^2=\sum\limits_{n\in\mathbb{N}^*}|a_n|^2\lambda_n^s<\infty\Big\},$$
 where $\dbN^*=\dbN\setminus\{0\}, $ and $\{\lambda_n\}_{n\in\dbN^*}$ and $\{\Phi_n\}_{n\in\dbN^*}$ are eigenvalues and eigenvectors of a degenerate elliptic operator, respectively.  Also, it is easy to check that $\mathcal{H}\subset L^2(0, 1)\times H_\alpha^*(0, 1),$ and   any initial value in the space $(L^2(0, 1)\times H_\alpha^*(0, 1))\setminus \mathcal{H}$ is not controllable in time $T$ (see \cite{Gueye}).

 Borrowing some ideas from  \cite{Gueye}, we study
the null controllability problems of one-dimensional degenerate
 wave equations. Different from \cite{Gueye},   the control acts on the nondegenerate boundary in our problems.
Notice that for the  one-dimensional nondegenerate wave equations,  the controllability results are same, no matter which endpoint of $[0, L]$ a control acts on. However, there is a  different new phenomena for degenerate wave equations. In fact,  we find that when a control enters the system from the nondegenerate boundary,   any initial value in the state space $L^2(0, L)\times H_\alpha^*(0, L)
$ is controllable.  But, when a control acted  on the degenerate boundary, only a subspace of the state space was  controllable in \cite{Gueye}.
Now, we explain it
from the physical background of degenerate wave equations.
Due to degeneracy, the propagation speed of waves is zero in the endpoint $x=0$. If a control acts on this point, it
has little effect on the whole waveform. But, when control acts on the nondegenerate boundary $x=L$,  it will give
a permanent  influence on the whole waveform. Therefore,  any  initial value in the state function  can be driven to zero.

\medskip

Furthermore,  by the method of  characteristic lines,  we explain   controllability times  for one-dimensional nondegenerate and degenerate wave equations, respectively.  Notice that a wave  propagates along  characteristic lines,  and for $\alpha\in [0, 2)$,
its  travel  time  from the endpoint $x=L$ to the other $x=0$ is calculated as follows:
\begin{equation}\label{zg2}
\int_0^L\frac{1}{x^{\frac{\alpha}{2}}}dx=\displaystyle\frac{T_\alpha}{2}.
\end{equation}
For the degenerate wave equation ($\alpha>0$), though the propagation speed of  the wave   tends to zero,  as it is sufficiently  close to the endpoint $x=0$,  its distance  to this point also  tends to $0$. Hence, the travel time  may be finite by (\ref{zg2}).
When a control acts  on the right endpoint $x=L$, it is easy to see that  the controllability time is $T_\alpha$ (resp. $+\infty$) for $\alpha\in [0, 2)$ (resp. $\alpha\geq 2$) (from Figures 1-3).

\begin{center}
\includegraphics[totalheight=60mm]{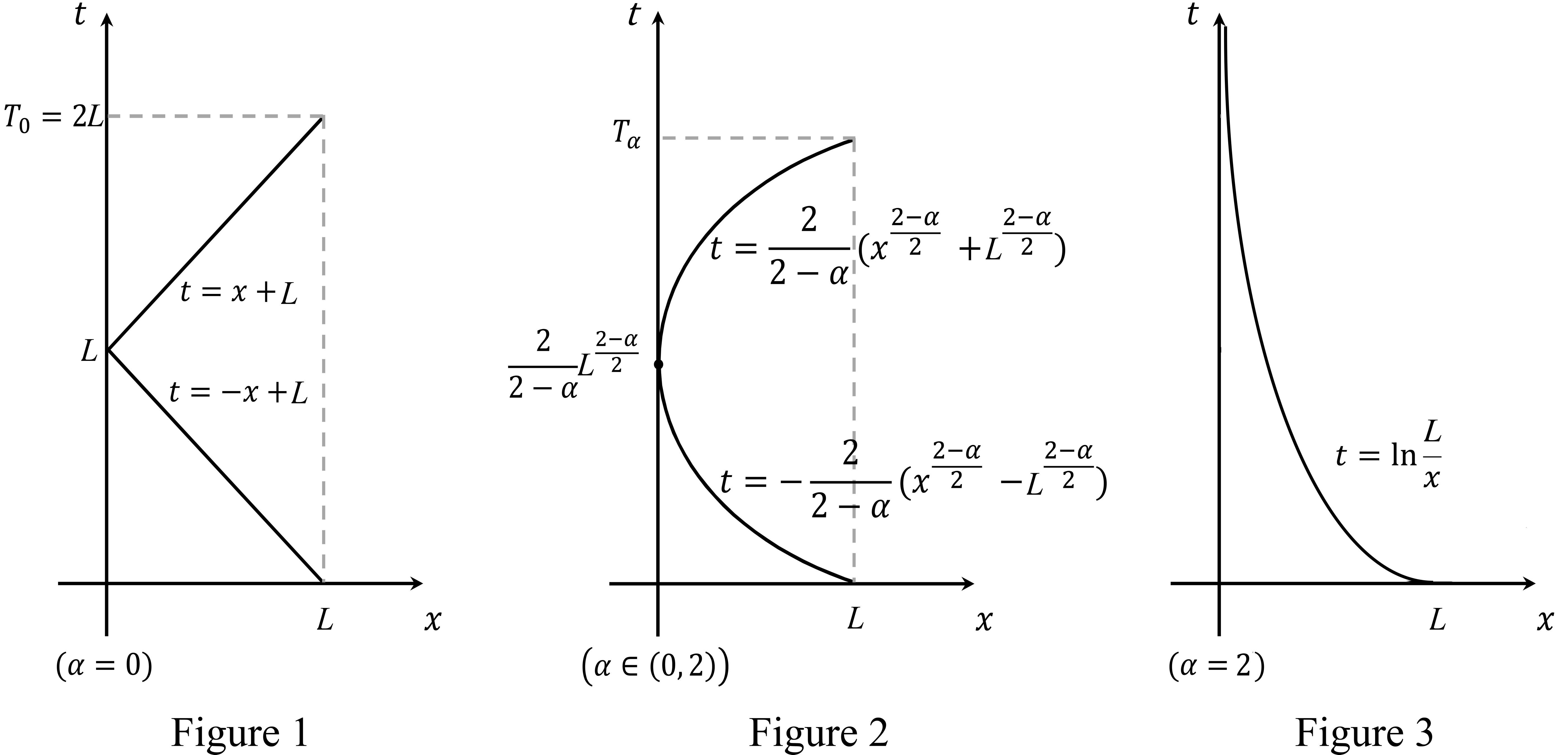}\ \
\end{center}

The rest of this paper is organized as follows. In Section 2,
 the well-posedness results for the system (\ref{1.1}) are given.
In Section 3, we prove the main controllability results of this paper (Theorem {\ref{t1}}).

\section{Well-posedness of degenerate wave equations}

In this section, we prove the well-posedness of the equation (\ref{1.1}) (Theorem \ref{t1.3}).
First,
consider the following linear degenerate wave equation:
\begin{equation}\label{1.3}
\left\{\begin{array}{ll}
y_{tt}-(x^\alpha y_x)_x=f & (x,t)\in Q,\\[2mm]
\left\{\begin{array}{ll}
y(0,t)=0 &(0<\alpha<1)\\[2mm]
(x^\alpha y_x)(0,t)=0 &(\alpha\geq1)\\[2mm]
\end{array}\right. &t\in(0,T),\\[2mm]
y(L,t)=0 &t\in(0,T),\\[2mm]
y(x,0)=y_0(x),\ y_t(x,0)=y_1(x) &x\in (0, L),
\end{array}\right.
\end{equation}
where $f\in L^1(0,T;L^2(0,L))$ and  $(y_0,y_1)\in H_\alpha^1(0,L)\times L^2(0,L)$.

\medskip

We are concerned with  weak solutions of the system $(\ref{1.3})$.

\begin{definition}\label{d1.1}
A function $y\in C([0,T];H_\alpha^1(0,L))\bigcap C^1([0,T];L^2(0,L))$ is said to be a weak solution of the system $(\ref{1.3})$, if
for any $\varphi\in L^{2}(0,T;H_\alpha^1(0,L))$ satisfying   $\varphi_t\in L^2(Q)$
and $\varphi(\cdot,T)=0$, it holds that $y(x, 0)=y_0(x)$ in $(0, L)$ and
\begin{eqnarray}\label{0.2}
\begin{array}{rl}
\displaystyle\int_{Q}(-y_t\varphi_t+x^\alpha y_x\varphi_x)dxdt-\displaystyle\int^L_0 y_1\varphi(x, 0)dx
=\displaystyle\int_{Q}f\varphi dxdt.\\[3mm]
\end{array}
\end{eqnarray}
\end{definition}

We have the following well-posedness result for the  equation (\ref{1.3}).
\begin{theorem}\label{t1.1}
For any $f\in L^1(0,T;L^2(0,L))$ and $(y_0,y_1)\in H_\alpha^1(0,L)\times L^2(0,L)$,
the system $(\ref{1.3})$ has a unique weak solution
$y\in C([0,T];H_\alpha^1(0,L))\cap C^1([0,T];L^2(0,L))$.
Moreover, there exists a constant $C=C(T,L,\alpha)$ such that
\begin{eqnarray*}
&&\|y\|_{L^\infty(0,T;L^2(0,L))}+\|y_t\|_{L^\infty(0,T;L^2(0,L))}+
\|x^\alpha y_x^2\|_{L^\infty(0,T;L^1(0,L))}\\[3mm]
&&\leq C\big{(}\|f\|_{L^1(0,T;L^2(0,L))}+\|y_0\|_{H^1_\alpha(0,L)}+
\|y_1\|_{L^2(0,L)}\big{)}.
\end{eqnarray*}
\end{theorem}

\noindent {\bf Proof. } We borrow some ideas from \cite{wang1}.
First,  suppose that  $f$ and $(y_0,y_1)$ are sufficiently smooth.
For any positive integer $n$, consider the following
nondegenerate wave equation:
\begin{equation}\label{1.4}
\left\{\begin{array}{ll}
y_{tt}-\(\big{(}x^\alpha+\frac{1}{n}\big{)}y_{x}\)_x=f, &(x,t)\in Q,\\[3mm]
y(0,t)=y(L,t)=0, &t\in(0,T),\\[3mm]
y(x,0)=y_0(x),\  y_t(x,0)=y_1(x), &x\in (0,L).
\end{array}\right.
\end{equation}
By the classical
theory of wave equations, the system (\ref{1.4}) admits a unique classical solution $y^n$.
 For simplicity,
we denote by $y$ the solution $y^n$.
Multiplying both sides of the first equation of (\ref{1.4})
by $y_t$ and integrating it in $Q_s=(0,s)\times (0,L)$,
by the H\"{o}lder inequality, we have that
\begin{eqnarray*}
\begin{array}{rl}
&\max\limits_{s\in[0,T]}\displaystyle\int_0^L\(y_t^2(x,s)+(x^\alpha+\frac{1}{n})|y_x(x,s)|^2\)dx\\[3mm]
&\leq \displaystyle\int_0^L\(y_t^2(x,0)+(x^\alpha+\frac{1}{n})|y_x(x,0)|^2\)dx
+\max\limits_{s\in[0,T]}\|y_t\|_{L^2(0,L)}\displaystyle\int_{0}^T\|f\|_{L^2(0,L)}dt.
\end{array}
\end{eqnarray*}
This implies that
\begin{eqnarray}\label{1.5}
\begin{array}{rl}
&\displaystyle\int_0^L\(y_t^2(x,s)+(x^\alpha+\frac{1}{n})|y_x(x,s)|^2\)dx\\[3mm]
&\leq C\displaystyle\int_0^L\(y_t^2(x,0)+(x^\alpha+\frac{1}{n})|y_x(x,0)|^2\)dx
+C\(\displaystyle\int_{0}^T\|f\|_{L^2(0,L)}dt\)^2.
\end{array}
\end{eqnarray}

On the other hand, notice that
$$y(x,s)=y(x,0)+\int_0^sy_t(x,t)dt.$$
By the H\"{o}lder inequality, we arrive at
\begin{equation}\label{1.6}
\displaystyle\int_0^Ly^2(x,s)dx
\leq 2\displaystyle\int_0^Ly^2(x,0)dx
+2s\displaystyle\int_{Q_s}y_t^2(x,t)dxdt.
\end{equation}
By (\ref{1.5}) and (\ref{1.6}), it follows that
\begin{eqnarray*}
\begin{array}{rl}
&\displaystyle\int_0^L\(y^2(x,s)+y_t^2(x,s)+(x^\alpha+\frac{1}{n})|y_x(x,s)|^2\)dx\\[3mm]
&\leq C\[\displaystyle\int_0^L\(y^2(x,0)+y_t^2(x,0)+(x^\alpha+\frac{1}{n})|y_x(x,0)|^2\)dx\\[3mm]
&\ \ \ +\displaystyle\int_{Q_s}\(y_t^2+y^2+(x^\alpha+\frac{1}{n})|y_x|^2\)dxdt+
\(\displaystyle\int_0^T\|f\|_{L^2(0,L)}dt\)^2\].
\end{array}
\end{eqnarray*}
By Gronwall's inequality, we get that
\begin{eqnarray}\label{1.7}
&&\|y\|_{L^\infty(0,T;L^2(0,L))}+\|y_t\|_{L^\infty(0,T;L^2(0,L))}+
\|(x^\alpha+\frac{1}{n})|y_x|^2\|_{L^\infty(0,T;L^1(0,L))}\nonumber\\[3mm]
&&\leq C\(\|f\|_{L^1(0,T;L^2(0,L))}+
\|y_0\|_{H^1(0,L)}+\|y_1\|_{L^2(0,L)}\).
\end{eqnarray}
Hence, there exist a subsequence $\{y^{n_j}\}$ of $\{y^n\}$, and
a function $y\in L^\infty(0,T;L^2(0,L))$ satisfying $y_t\in L^\infty(0,T; L^2(0,L))$, such that as $j\rightarrow \infty$,
\begin{eqnarray}\label{zg00}
\begin{array}{rl}
&y^{n_j}\rightarrow y\ \ \ \mbox{weakly}\mbox{ in } L^2(Q);\ \ \ y_t^{n_j}\rightarrow y_t\ \ \ \mbox{weakly}\text{ in } L^2(Q);\\[4mm]
& \sqrt{x^\alpha+\frac{1}{n_j}}y_x^{n_j} \rightarrow \sqrt{x^\alpha}y_x\ \ \ \mbox{weakly}\text{ in } L^2(Q).
\end{array}
\end{eqnarray}
Since $\{y^{n_j}\}$ is the classical solution of (\ref{1.4}),
   we have that $y^{n_j}(x, 0)=y_0(x)$, and for any $\varphi\in C^\infty(\overline{Q})$, with $\varphi=0$ in some neighborhood of $\{0\}\times(0, T)$, $\{L\}\times(0, T)$  and $(0, L)\times\{T\}$,
\begin{eqnarray*}
\begin{array}{rl}
&\displaystyle\int_{Q}\Big[-y^{n_j}_t\varphi_t+(x^\alpha+\frac{1}{n_j}) y_x\varphi_x\Big]dxdt-\displaystyle\int_0^Ly_1(x)\varphi(x,0)dx
=\displaystyle\int_{Q}f\varphi dxdt.\\[3mm]
\end{array}
\end{eqnarray*}
By (\ref{zg00}), taking a limit in the above equality and noticing that $C^\infty_0(Q)$ is dense in $L^2(0, T; H^1_\alpha(0, L))$, we obtain that $y$ satisfies (\ref{0.2}) for any $\varphi\in L^2(0, T; H^1_\alpha(0, L))$ with  $\varphi_t\in L^2(Q)$
and $\varphi(\cdot,T)=0$.
Also,  it is easy to show that $\{y^{n_j}\}$ is a Cauchy sequence in $C([0,T];H_\alpha^1(0,L))$ $\bigcap C^1([0,T];L^2(0,L))$. This implies
$y\in  C([0,T];H_\alpha^1(0,L))\bigcap C^1([0,T];L^2(0,L))$.
Hence,  $y$ is the weak solution of the system (\ref{1.3}) associated to  smooth functions $(y_0,y_1)$ and $f$.

\medskip
Next, we prove the existence of weak solutions of (\ref{1.3}) for any $(y_0,y_1)\in H_\alpha^1(0,L)\times L^2(0,L)$ and $f\in L^1(0, T; L^2(0, L))$. Let $\{y_0^m\}$, $\{y_1^m\}$ and $\{f^m\}$ be sequences of smooth functions, respectively,
such that as $m\rightarrow \infty$,
\begin{equation}\label{ff}
y_0^m\rightarrow y_0\text{  in } H_\alpha^1(0,L),\quad  y_1^m\rightarrow y_1 \text{ in } L^2(0,L)\quad\mbox{ and }\quad
f^m\rightarrow f\mbox{ in } L^1(0, T; L^2(0, L)).
\end{equation}
Denote by $y^m$ the solution of (\ref{1.3}) associated to  $(y_0^m,y_1^m)$ and $f^m$.  Similarly, we can obtain that
\begin{eqnarray}\label{e}
&&\|y^m-y^n\|_{C([0,T];H_\alpha^1(0,L))}+\|y^m_t-y^n_t\|_{C([0,T];L^2(0,L))}\nonumber\\[3mm]
&&\leq C\(\|f^m-f^n\|_{L^1(0,T;L^2(0,L))}+
\|y_0^m-y_0^n\|_{H^1_\alpha(0,L)}+\|y_1^m-y_1^n\|_{L^2(0,L)}\).
\end{eqnarray}
Therefore,  there exists  $y\in C([0,T];H_\alpha^1(0,L))\bigcap C^1([0,T];L^2(0,L))$, such that as    $m\rightarrow\infty$,
$$y^m\rightarrow y\ \ \ \text{in}\ \ \ C([0,T];H_\alpha^1(0,L))\ \ \text{and}\  \
y^m_t\rightarrow y_t\ \ \ \text{in}\ \ \ C([0,T];L^2(0,L)).$$
Moreover, it is easy to show that  $y$ is the weak solution of (\ref{1.3}) for $(y_0,y_1)\in H_\alpha^1(0,L)\times L^2(0,L)$ and $f\in L^1(0, T; L^2(0, L))$.
\medskip

Finally, we prove
the uniqueness of weak solutions. Let $\widetilde{y}$ and $\overline{y}$ be two weak solutions of the system (\ref{1.3}), and set
 $\widehat{y}=\widetilde{y}-\overline{y}$. Then
$\widehat{y}\in C([0,T];H_\alpha^1(0,L))\cap C^1([0,T];L^2(0,L))$ satisfies
\begin{eqnarray}\label{1.8}
\displaystyle\int_{Q}(-\widehat{y}_t\varphi_t+x^\alpha\widehat{y}_x\varphi_x)
dxdt=0,
\ \forall\ \varphi\in L^{2}(0,T;H_\alpha^1(0,L))\mbox{ with } \varphi_t\in L^2(Q) \mbox{ and }
\varphi(\cdot,T)=0.
\end{eqnarray}
For any $\ell\in C_0^\infty(Q)$, consider the following degenerate wave equation:
\begin{eqnarray}\label{zg3}
\left\{\begin{array}{ll}
\psi_{tt}-(x^\alpha\psi_x)_x=\ell &(x,t)\in Q,\\[3mm]
\left\{\begin{array}{ll}
\psi(0,t)=0 &(0<\alpha<1)\\[2mm]
(x^\alpha\psi_x)(0,t)=0&(\alpha\geq1)\\[2mm]
\end{array}\right.&t\in(0,T),\\[2mm]
\psi(L,t)=0&t\in(0,T),\\[2mm]
\psi(x,T)=0,\ \psi_t(x,T)=0 &x\in(0,L).
\end{array}\right.
\end{eqnarray}
By the above argument on the existence of weak solutions,  one can find a weak solution $\psi\in C([0,T];H_\alpha^1(0,L))\cap C^1([0,T];L^2(0,L))$ of (\ref{zg3}). This implies that
\begin{eqnarray}\label{11}
&\displaystyle\int_{Q}(-\psi_t\varphi_t+x^\alpha\psi_x\varphi_x)
dxdt=\displaystyle\int_{Q}\ell\varphi dxdt,\\
&\ \forall\ \varphi\in L^{2}(0,T;H_\alpha^1(0,L))\mbox{ with }  \varphi_t\in L^2(Q)\mbox{ and }
 \varphi(\cdot,0)=0\ \text{in}\ (0,L).\nonumber
\end{eqnarray}
Choosing $\varphi=\psi$ in (\ref{1.8}),
and $\varphi=\widehat{y}$ in (\ref{11}),
  we have that
\begin{equation*}
\displaystyle\int_{Q}\ell\widehat{y}dxdt=0,\ \forall\ \ell\in C_0^\infty(Q).
\end{equation*}
Hence,  $\widehat{y}(x,t)=0$ a.e. in $Q$. The proof is completed.\endpf

\medskip

\medskip

\begin{remark}
The
 well-posedness of the system $(\ref{1.3})$ was proved  by the semigroup theory $($for
$\alpha\in(0,1))$ in $\cite[Proposition\ 4.2]{Gueye}$.
Here we use another method to prove it for any $\alpha\in(0,\infty)$.  Our method is also applicable to  the more general multi-dimensional degenerate wave operator:
$$\mathcal{L}y=y_{tt}-\mbox{div}(b(x, t) \nabla y),  $$
where $b\in C(\overline{\Omega}\times [0,T])$ satisfies $\displaystyle\frac{b_t}{b}\in L^\infty(\Omega\times (0,T))$ and
$b>0$ in $\Omega\times(0,T)$, with $\Omega$ being a nonempty bounded domain of $\mathbb{R}^n$.

\end{remark}

In the following, we prove
a hidden regularity result for solutions of the following
degenerate wave equation:
\begin{equation}\label{1.11}
\left\{\begin{array}{ll}
z_{tt}-(x^\alpha z_x)_x=h &(x,t)\in Q,\\[2mm]
\left\{\begin{array}{ll}
z(0,t)=0 &(0<\alpha<1)\\[2mm]
(x^\alpha z_x)(0,t)=0 &(\alpha\geq1)\\[2mm]
\end{array}\right.&t\in(0,T),\\[2mm]
z(L,t)=0&t\in(0,T),\\[2mm]
z(x,0)=z_0(x),\ z_t(x,0)=z_1(x) &x\in (0,L),
\end{array}\right.
\end{equation}
where $h\in L^1(0,T;L^2(0,L))$ and $(z_0,z_1)\in H_\alpha^1(0,L)\times L^2(0,L).$ Define the following energy functional:
$$E(t)=\frac{1}{2}\int_0^L\[z_t(x,t)^2+x^\alpha z_x^2(x, t)\]dx.$$
As a preliminary,  we have the following energy estimate.
\begin{lemma}\label{l1.1}
There exists a constant $C>0$ such that
$E(t)\leq C\( E(0)+\Big(\displaystyle\int_0^T\|h\|_{L^2(0,L)}dt\Big)^2
\),$ $\forall\ t\in[0,T].$
\end{lemma}

 \noindent {\bf Sketch of the proof. }
  For any $h\in L^1(0,T;L^2(0,L))$,
multiplying both sides of the first equation in (\ref{1.11}) by
$z_t$ and integrating it on $(0,L)$, by Gronwall's inequality and a simple calculation,
we can obtain the desired estimate. \endpf

\medskip

By Lemma \ref{l1.1}, we have the following hidden regularity result for
 (\ref{1.11}).
\begin{proposition}\label{t1.2}
For any $\alpha\in(0,2)$, any solution $z$ of $(\ref{1.11})$ satisfies $z_x(L,\cdot)\in L^2(0,T).$
Moreover,
\begin{eqnarray}\label{1.12}
\int_0^Tz_x^2(L,t)dt\leq C\(E(0)+\Big(\displaystyle\int_0^T\|h\|_{L^2(0,L)}dt\Big)^2
\).
\end{eqnarray}
\end{proposition}

\noindent {\bf Proof. } It suffices to prove (\ref{1.12}) for classical solutions  of (\ref{1.11}).
Hence,  we assume that $h$ and $(z_0,z_1)$ are sufficiently smooth.  Then the system (\ref{1.11}) admits a unique classical solution $z\in C([0,T];H_\alpha^2(0,L))\cap C^1([0,T];H_\alpha^1(0,L))$ (see \cite[Proposition\ 4.2]{Gueye}), where
$H_\alpha^2(0,L)=\Big\{u \in H_\alpha^1(0,L)\big|\ x^\alpha u_x\in H^1(0,L)\Big\}$.
Choose $
q(x)=
\left\{\begin{array}{ll}
x &x\in[0,\frac{L}{2}),\\
-2x+\frac{3}{2}L&x\in[\frac{L}{2}, L],
\end{array}\right.$
multiplying  both sides of the first equation in (\ref{1.11}) by
$q z_x$ and integrating it on $Q$,  we obtain that
\begin{eqnarray}\label{zg1}
\begin{array}{rl}
\displaystyle\frac{1}{2}\int_0^Tx^{\alpha}q z_x^2dt\Big|_0^L
=&\displaystyle\int_0^L qz_tz_xdx\Big|_0^T+
\frac{1}{2}\int_Qq_x\big{(}z_t^2+x^\alpha z_x^2\big{)}dxdt
\displaystyle-\frac{1}{2}\int_Q\alpha qx^{\alpha-1}z_x^2dxdt\\[3mm]
&-\displaystyle\int_Qqhz_xdxdt-\displaystyle\frac{1}{2}\int_0^Tqz_t^2dt\Big|_0^L.
\end{array}
\end{eqnarray}
(\ref{zg1}) implies that
{\setlength\arraycolsep{2pt}
\begin{eqnarray}\label{1.13}
\begin{array}{rl}
&\displaystyle\frac{L^{\alpha+1}}{2}\int_0^Tz_x^2(L,t)dt+\int_0^T(x^{\alpha+1}z_x^2)(0,t)dt\\
&
=-2\displaystyle\int_0^L qz_tz_xdx\Big{|}_0^T-
\int_Qq_x\(z_t^2+x^\alpha z_x^2\)dxdt+\displaystyle\int_Q\alpha qx^{\alpha-1}z_x^2dxdt\\
&\quad
+2\displaystyle\int_Qqhz_xdxdt
+\displaystyle\int_0^Tqz_t^2dt\Big|_0^L.
\end{array}
\end{eqnarray}}
When $\alpha\in(0,2)$, $(xz^2)(0,t)=0$ (see \cite[Proposition\ 2.4]{car}). Therefore,
$
\displaystyle\int_0^Tqz_t^2dt\Big|_0^L=0.
$

Further,
by the H\"{o}lder inequality and Lemma \ref{l1.1}, we arrive at
{\setlength\arraycolsep{2pt}
\begin{eqnarray}\label{1.14}
&&-2\displaystyle\int_0^L qz_tz_xdx\Big{|}_0^T
=-2\displaystyle\int_0^{\frac{L}{2}} xz_tz_xdx\Big{|}_0^T
-2\displaystyle\int_{\frac{L}{2}}^L z_t(-2x+\frac{3}{2}L)z_xdx\Big{|}_0^T\nonumber\\
&&\leq C\(E(0)+\(\displaystyle\int_0^T\|h\|_{L^2(0,L)}dt\)^2
\).
\end{eqnarray}}
Similarly,
{\setlength\arraycolsep{2pt}
\begin{eqnarray}\label{1.15}
&&-\displaystyle\int_Qq_x\(z_t^2+x^\alpha z_x^2\)dxdt+\displaystyle\int_Q\alpha q x^{\alpha-1}z_x^2dxdt
+2\displaystyle\int_Qqhz_xdxdt\nonumber\\
&&\leq C\(E(0)+\(\displaystyle\int_0^T\|h\|_{L^2(0,L)}dt\)^2
\).
\end{eqnarray}}
Combining  (\ref{1.14}) and (\ref{1.15}) with (\ref{1.13}),
we get the desired estimate (\ref{1.12}).\endpf

\medskip

\medskip

Now, we give a proof of Theorem \ref{t1.3}.

\medskip

\noindent {\bf Proof of Theorem \ref{t1.3}. }
First, by Theorem \ref{t1.1} and Proposition \ref{t1.2}, the equation
(\ref{1.9}) admits a unique weak solution $v\in C([0,T];H_\alpha^1(0,L))\cap C^1([0,T];L^2(0,L))$
satisfying $v_x(L,\cdot)\in L^2(0,T).$

Next, define a linear functional on $L^1(0,T;L^2(0,L))$:
$$\mathscr{L}(\xi)=\langle w_1,v(0)\rangle _{H_\alpha^*(0,L),H_\alpha^1(0,L)}
-\int_0^L w_0v_t(0)dx-L^\alpha\int_0^T\theta(t)v_x(L,t)dt,\ \ \forall\ \xi\in L^1(0,T;L^2(0,L)).$$
Then we have
$$|\mathscr{L}(\xi)|\leq|(w_0,w_1)|_{L^2(0,L)\times H_\alpha^*(0,L)}
\cdot|(v(0),v_t(0))|_{H_\alpha^1(0,L)\times L^2(0,L)}
+C\|\theta\|_{L^2(0,T)}\cdot\|v_x(L,\cdot)\|_{L^2(0,T)}.$$
By Proposition \ref{t1.2} and Lemma \ref{l1.1}, $\mathscr{L}$ is a bounded linear functional on $L^1(0,T;L^2(0,L))$.
Therefore, there exists a  function $w\in L^\infty(0,T;L^2(0,L))$, such that  (\ref{1.10}) holds.
By the standard smoothing technique, we can get $w\in C([0,T];L^2(0,L))$.

Similarly, consider the following degenerate wave equation:
\begin{equation*}
\left\{\begin{array}{ll}
v_{tt}-(x^\alpha v_x)_x=g_t& (x,t)\in Q,\\[2mm]
\left\{\begin{array}{ll}
v(0,t)=0&(0<\alpha<1)\\[2mm]
(x^\alpha v_x)(0,t)=0 &(\alpha\geq1)\\[2mm]
\end{array}\right.&t\in(0,T),\\[2mm]
v(L,t)=0&t\in(0,T),\\[2mm]
v(x,T)=0,\ v_t(x,T)=0&x\in (0,L),
\end{array}\right.
\end{equation*}
where $g\in C_0^\infty(0,T; H_\alpha^1(0,L))$. For this system, we can
obtain the conclusions similar to  Proposition \ref{t1.2} and Lemma \ref{l1.1}. Notice that the density of $C_0^\infty(0,T; H_\alpha^1(0,L))$ in $L^1(0,T; H_\alpha^1(0,L))$,
we can also prove $w\in C^1([0,T];H_\alpha^*(0,L))$.
This completes the proof.\endpf

\section{Controllability of degenerate wave equations}

In this section, we study the boundary controllability of the
degenerate wave equation (\ref{1.1}).
First,  consider the following linear degenerate wave equation:

\begin{equation}\label{2.1}
\left\{\begin{array}{ll}
v_{tt}-(x^\alpha v_x)_x=0 & (x,t)\in Q,\\[2mm]
\left\{\begin{array}{ll}
v(0,t)=0 &(0<\alpha<1)\\[2mm]
(x^\alpha v_x)(0,t)=0 &(\alpha\geq1)\\[2mm]
\end{array}\right.&t\in(0,T),\\[2mm]
v(L,t)=0&t\in(0,T),\\[2mm]
v(x,0)=v_0(x),\ v_t(x,0)=v_1(x)&x\in (0,L),
\end{array}\right.
\end{equation}
where $(v_0,v_1)\in H_\alpha^1(0,L)\times L^2(0,L)$.

By the  duality technique, it is easily seen that, the null controllability  of (\ref{1.1})
can be reduced to an observability estimate for the system
(\ref{2.1}).
\begin{proposition}\label{p1}
The system $(\ref{1.1})$ is null controllable in time
$T$ with a control $\theta\in L^2(0, T)$, if and only if there exists a constant $C>0$, such that any solution $v$ of $(\ref{2.1})$ satisfies
\begin{equation}\label{2.2}
\|v_0\|^2_{H_\alpha^1(0,L)}+\|v_1\|^2_{L^2(0,L)}
\leq C\int_0^Tv_x^2(L,t)dt,\ \ \forall \ (v_0,v_1)\in H_\alpha^1(0,L)\times L^2(0,L).
\end{equation}
\end{proposition}

\medskip

First,  introduce some notations. We write $A\asymp B$,  if there exist two constants $C_1, C_2>0$,
such that $C_1A\leq B\leq C_2A$.
 Set
$ l^2(\mathbb{Z})=\Big\{\{a_n\}_{n\in\mathbb{Z}} \ | \ a_n\in \mathbb{C},\
\sum\limits_{n\in\mathbb{Z}}|a_n|^2<+\infty \Big\}$ and $i$ denotes the imaginary unit.

In order to prove  the observability inequality (\ref{2.2}) for the system (\ref{2.1}),  we adopt the similar method  used
 in \cite[Theorem\ 3.1]{Gueye}, based on the following  known Ingham inequality.

\begin{lemma}\label{l1.2}$(\cite[Theorem\ 9.2]{kom})$
Suppose that $\{\eta_n\}_{n\in\mathbb{Z}}$ is a sequence of real numbers satisfying the following gap condition  for some  constant $\delta>0$,
$$\eta_{n+1}-\eta_n\geq\delta,\ \ \forall\  n\in \mathbb{Z}.$$ Let
 $n^+(r)$ denote the largest number of terms of the sequence
$\{\eta_n\}_{n\in\mathbb{Z}}$ contained in an interval of length $r$, and $D^+=\lim\limits_{r\rightarrow\infty}\frac{n^+(r)}{r}$.
Then, it holds that

(1) For any $T>2\pi D^+$,
\begin{equation}\label{2.5}
\int_0^T\big{|}\sum\limits_{n\in\mathbb{Z}}a_ne^{i\eta_nt}\big{|}^2dt\asymp
\sum\limits_{n\in\mathbb{Z}}\big{|}a_n\big{|}^2, \quad
\mbox{for any } \{a_n\}_{n\in\mathbb{Z}}\in l^2(\mathbb{Z}).
\end{equation}

(2) For any  $T<2\pi D^+$, $(\ref{2.5})$ does not hold.
\end{lemma}

The main result of this section is the following observability inequality for the system (\ref{2.1}).

\begin{theorem}\label{t2.1}
Suppose that $\alpha\in (0,2)$. Then for any $T>T_\alpha$ $($defined in $(\ref{zg}))$,
 the estimate $(\ref{2.2})$ holds for  $(\ref{2.1})$.  For any $T<T_\alpha$,
 the estimate $(\ref{2.2})$ does not hold.
\end{theorem}
\noindent {\bf Proof.}
(1) For any
$\alpha\in(0,1),$ we
consider the eigenvalue problem:
\begin{equation}\label{2.10}
\left\{\begin{array}{ll}
-(x^{\alpha}\Phi'(x))'=\lambda \Phi(x)\ \ x\in(0,L),\\[2mm]
\Phi(0)=\Phi(L)=0.\\[2mm]
\end{array}\right.
\end{equation}
Similar to
$\cite{Gueye}$,   for any $n\geq 1,$  the solutions of (\ref{2.10}) are given as follows:
\begin{equation}\label{2.15}
\lambda_n=\Big(\frac{\rho j_{\mu,n}}{L^\rho}\Big)^2 \quad\mbox{and}
\quad
\Phi_n(x)=\frac{(2\rho)^{\frac{1}{2}}}
{L^\rho|J_{\mu}^{'}(j_{\mu,n})|}x^{\frac{1-\alpha}{2}}J_{\mu}\Big(j_{\mu,n}\frac{x^\rho}{L^\rho}\Big),
\ \forall\  x\in (0,L),
\end{equation}
where
\begin{equation}\label{2.17}
\mu=\displaystyle\frac{1-\alpha}{2-\alpha},\ \rho=\displaystyle\frac{2-\alpha}{2},\
J_{\mu}(x)=\sum\limits_{m\geq0}\frac{(-1)^m}{m!\cdot\Gamma(m+\mu+1)}
\Big(\frac{x}{2}\Big)^{2m+\mu},\  x\geq0,
\end{equation}
$\Gamma(\cdot)$ is the Gamma function,
and $\{j_{\mu,n}\}_{n\geq1}$ are the positive zeros of the Bessel function
$J_{\mu}$.
Moreover, for any $\mu\geq-\displaystyle\frac{1}{2}$ and $\beta>0$,
\begin{equation}\label{2.6}
\int_0^Lx^{2\beta-1}J_{\mu}\Big(j_{\mu,n}(\frac{x}{L})^\beta\Big)J_{\mu}\Big(j_{\mu,m}(\frac{x}{L})^\beta\Big)dx=
L^{2\beta}\frac{\delta_{n}^m}{2\beta}[J'_{\mu}(j_{\mu,n})]^2,
\end{equation}
where $\delta_{n}^m$ is the Kronecker symbol (see \cite[(4.26)]{Gueye}).
 Therefore, $\{\Phi_n\}_{n\in\mathbb{N}^*}$ is an orthonormal basis of $L^2(0,L)$.

\medskip

Next, for any $(v_0,v_1)\in H_\alpha^1(0,L)\times L^2(0,L)$,  set
\begin{equation*}
v_0(x)=\sum\limits_{n\in\mathbb{N}^*}v_0^n\Phi_n(x)\quad\mbox{ and }\quad
v_1(x)=\sum\limits_{n\in\mathbb{N}^*}v_1^n\Phi_n(x).
\end{equation*}
Then it is easy to check that the solution $v$ of (\ref{2.1})
 is as follows:
\begin{equation*}
v(x,t)=\sum\limits_{n\in\mathbb{N}^*}v_n(t)\Phi_n(x)\quad\mbox{with }
v_n(t)=c_ne^{i\sqrt{\lambda_n}t}+c_{-n}e^{-i\sqrt{\lambda_n}t},
\end{equation*}
where 
\begin{equation*}
c_n=\frac{1}{2}\Big(v_0^n+\frac{v_1^n}{i\sqrt{\lambda_n}}\Big)\quad \mbox{ and }\quad
c_{-n}=\frac{1}{2}\Big(v_0^n-\frac{v_1^n}{i\sqrt{\lambda_n}}\Big).
\end{equation*}
By a simple calculation, we have that
\begin{equation*}
\Phi'_n(x)=\frac{(2\rho)^{\frac{1}{2}}}{L^\rho|J_{\mu}^{'}(j_{\mu,n})|}
\Big(\frac{1-\alpha}{2}x^{\frac{-1-\alpha}{2}}J_{\mu}(j_{\mu,n}\frac{x^\rho}{L^\rho})
+x^{\frac{1-2\alpha}{2}}J'_{\mu}(j_{\mu,n}\frac{x^\rho}{L^\rho})\frac{j_{\mu,n}\rho}{L^\rho}\Big).
\end{equation*}
It follows that
\begin{equation*}
\Phi'_n(L)=L^{-\frac{3}{2}}\frac{\sqrt{2\rho}}{|J_{\mu}^{'}(j_{\mu,n})|}
\rho j_{\mu,n}J_{\mu}^{'}(j_{\mu,n}).
\end{equation*}
This implies that
\begin{equation*}
v_x(L,t)=L^{-\frac{3}{2}}\rho\sqrt{2\rho}\sum\limits_{n\in\mathbb{N}^*}
j_{\mu,n}\frac{J_{\mu}^{'}(j_{\mu,n})}{|J_{\mu}^{'}(j_{\mu,n})|}
(c_ne^{i\sqrt{\lambda_n}t}+c_{-n}e^{-i\sqrt{\lambda_n}t}).
\end{equation*}
By Lemma \ref{l1.2}, we have that for any $T>2\pi D^+$,
\begin{eqnarray}\label{2.16}
&&\int_0^Tv_x^2(L,t)dt\asymp \sum\limits_{n\in\mathbb{N}^*}(j_{\mu,n})^2
\Big(|v_0^n|^2+\frac{|v_1^n|^2}{(\rho j_{\mu,n})^2}\Big)\nonumber\\[2mm]
&&\asymp \sum\limits_{n\in\mathbb{N}^*}(\lambda_n|v_0^n|^2+|v_1^n|^2)
\asymp \|v_0\|^2_{H_\alpha^1(0,L)}+\|v_1\|^2_{L^2(0,L)}.
\end{eqnarray}
Moreover, by the definition of $n^+(r)$,
$$\lim\limits_{r\rightarrow\infty}\[\frac{r}{n^+(r)}-
\frac{1}{n^+(r)}\sum\limits_{n=1}^{n^+(r)-1}\frac{\rho}{L^\rho}\int_{j_{\mu,n}}^{j_{\mu,n+1}}dx\]=0.$$
Since $(j_{\mu,n+1}-j_{\mu,n})$
converges to $\pi$ as $n\rightarrow +\infty$ (\cite[Lemma\ 1]{Gueye}),  we get that $D^+=\displaystyle\frac{L^\rho}{\rho\pi}.$ Hence,
(\ref{2.2}) holds for any
$T>2\pi D^+=T_\alpha$, but it fails  for any  $T<T_\alpha$.

\medskip

\medskip

\noindent (2) For any $\alpha\in[1,2),$ we
consider the eigenvalue problem:
\begin{equation}\label{2.21}
\left\{\begin{array}{ll}
-(x^{\alpha}\widehat{\Phi'}(x))'=\lambda \widehat{\Phi}(x)\ \ x\in(0,L),\\[2mm]
(x^{\alpha}\widehat{\Phi}')(0)=\widehat{\Phi}(L)=0.\\[2mm]
\end{array}\right.
\end{equation}
 Similar to the case of $\alpha\in(0, 1)$,  the  solutions  of (\ref{2.21})
are as follows:
\begin{equation*}
\widehat{\lambda}_n=\Big(\frac{\rho j_{\widehat{\mu},n}}{L^\rho}\Big)^2,\ n\geq1\quad\mbox{and}\quad
\widehat{\Phi}_n(x)=\frac{(2\rho)^{\frac{1}{2}}}
{L^\rho|J_{\widehat{\mu}}^{'}(j_{\widehat{\mu},n})|}x^{\frac{1-\alpha}{2}}J_{\widehat{\mu}}
\Big(j_{\widehat{\mu},n}\frac{x^\rho}{L^\rho}\Big),
\ \forall\ x\in (0,L),
\end{equation*}
where $\widehat{\mu}=\displaystyle\frac{\alpha-1}{2-\alpha}$ and $\rho=\displaystyle\frac{2-\alpha}{2}$.
Also, it is easy to check that  $\{\widehat{\Phi}_n\}_{n\in\mathbb{N}^*}$ is an orthonormal basis of $L^2(0,L)$.
The remainder  of the  proof  is similar to that in the case (1) and we omit it here. \endpf

\begin{remark}
It is well known that the controllability time $T=2L$ is optimal
for the classical one-dimensional wave equation $($with constant coefficients$)$ in $(0, L)$.
However,   the controllability of the system $(\ref{1.1})$ in time $T=T_\alpha$ is still an open problem.\end{remark}

\noindent {\bf Proof of (2) in Theorem \ref{t1}.}
 For any $\alpha>2$, let
$$X=\displaystyle\int_x^Ly^{-\frac{\alpha}{2}}dy\quad \text{and}\quad  W(X,t)=x^{\frac{\alpha}{4}}w(x,t)\quad  (\mbox{see }\cite{rd}).$$
Then the null  controllability of  the equation (\ref{1.1}) is transformed into a null controllability problem for the following nondegenerate wave equation
on the half-line:
\begin{equation}\label{c}
\left\{\begin{array}{ll}
W_{tt}(X,t)\!-\!W_{XX}(X,t)\!+\!M(X)W(X,t)\!=\!0& (X,t)\!\in\! \mathbb{R}^+\!\times\!(0,T),\\[2mm]
W(0,t)=\theta(t)&t\in(0,T),
\end{array}\right.
\end{equation}
where $M(X)=\frac{\alpha(3\alpha-4)}{\big[4L^{\frac{2-\alpha}{2}}-2(2-\alpha)X\big]^2}$.

On the other hand,
 by the duality, the null controllability of $(\ref{c})$ is equivalent to the following
observability estimate:
\begin{equation}\label{f}
|\psi_0|^2_{H_0^1(\mathbb{R}^+)}+
|\psi_1|^2_{L^2(\mathbb{R}^+)}
\leq C\int_0^T\psi_X^2(0,t)dt,
\end{equation}
where $\psi$ is the   solution of the nondegenerate wave equation:
\begin{equation}\label{l}
\left\{\begin{array}{ll}
\psi_{tt}-\psi_{XX}+M(X)\psi=0& (X,t)\in \mathbb{R}^+\times(0,T),\\[2mm]
\psi(0,t)=0&t\in(0,T),\\[2mm]
\psi(X,0)=\psi_0(X),\ \psi_t(X,0)=\psi_1(X)&X\in \mathbb{R}^+.
\end{array}\right.
\end{equation}

For  any given $\widehat{\psi}_0,
\widehat{\psi}_1\in C_0^\infty(\mathbb{R}^+)$, we choose
$\psi_0^n(X)=\widehat{\psi}_0(X-n)$ and  $\psi_1^n(X)=\widehat{\psi}_1(X-n)$ with sufficiently large  $n$ (see \cite{zua}).
Let $\psi^n$ be the solution of $(\ref{l})$ with initial value
$(\psi_0^n,\psi_1^n).$
\vspace{4mm}
It is easy to see that
\begin{equation*}
\frac{|\psi^n_0|^2_{H_0^1(\mathbb{R}^+)}+
|\psi_1^n|^2_{L^2(\mathbb{R}^+)}}{\displaystyle\int_0^T\psi_X^n(0,t)^2dt}\rightarrow\infty,
\ \ \text{as}\ n\rightarrow\infty.
\end{equation*}
This implies that the system $(\ref{c})$ is not null controllable.

\medskip

For $\alpha=2$,   we can also prove a similar result  for $M(X)=\frac{1}{4}$.
\endpf

\medskip

\medskip

$\mathbf{Note. }$
After completion of this work, we learned that Professor F. Alabau-Boussouira, Professor P. Cannarsa and Professor G. Leugering  proved the  null controllability for some one-dimensional degenerate wave equations with a more general diffusion coefficient, using the multiplier method in  the paper [F. Alabau-Boussouira, P. Cannarsa, G. Leugering, {\it Control and stabilization of degenerate wave equations},  arXiv:1505.05720]. Parts results of our work were contained in this paper.  But we prove the null controllability results of degenerate wave equations by a (different)  spectral method. Also,  the  controllability time $T_\alpha$ is sharp  in our work, i.e., if $T>T_\alpha$, the controllability holds and it is false if $T<T_\alpha$.
 Our work was reported by Professor Hang Gao  on  the 9th Workshop on Control of Distributed Parameter Systems on July 2, 2015 in Beijing.

\end{document}